\documentclass[12pt, a4paper]{amsart}

\usepackage{amscd,amssymb,amsthm}
\usepackage{graphicx}
\usepackage{color}
\usepackage{tikz,pgfplots}
\usepackage[colorlinks=true, pdfstartview=FitV, linkcolor=blue, citecolor=blue, urlcolor=blue]{hyperref}
\usepackage{wrapfig}
\usepackage{quoting}
\usepackage{stmaryrd}
\usepackage{hyperref}
\usepackage{multicol}
\textwidth \calclayout
\quotingsetup{font=small}

\definecolor{mypink1}{rgb}{0.858, 0.188, 0.478}
\definecolor{Prima}{rgb}{0.90,1,1}
\definecolor{Primab}{rgb}{0.96,1,1}
\definecolor{Seconda}{rgb}{1,0.90,1}
\definecolor{Segondab}{rgb}{1,0.96,1}
\definecolor{Terza}{rgb}{1,1,0.90}
\definecolor{Terzab}{rgb}{1,1,0.96}
\definecolor{Quarta}{rgb}{0.90,1,0.90}
\definecolor{Quartab}{rgb}{0.96,1,0.96}
\definecolor{yellow}{rgb}{1,1,0}
\definecolor{viola}{rgb}{1,0,1}
\definecolor{azzurro}{rgb}{0,1,1}
\definecolor{verde}{rgb}{0,1,0}
\definecolor{rosso}{rgb}{1,0,0}
\definecolor{blue}{rgb}{0,0,1}
\definecolor{ddd}{rgb}{1,1,0.1}

\newcommand{\beq}{\begin{equation}}
\newcommand{\eeq}{\end{equation}}

\theoremstyle{plain}
\newtheorem{theorem}{Theorem}

\theoremstyle{definition}

\theoremstyle{remark}
\newtheorem{remark}[theorem]{Remark}

\numberwithin{equation}{section}

\numberwithin{table}{section}

\definecolor{Primera}{rgb}{0.90,1,1}
\definecolor{Primerb}{rgb}{0.96,1,1}
\definecolor{Segona}{rgb}{1,0.90,1}
\definecolor{Segonb}{rgb}{1,0.96,1}
\definecolor{Tercera}{rgb}{1,1,0.90}
\definecolor{Tercerb}{rgb}{1,1,0.96}
\definecolor{Quarta}{rgb}{0.90,1,0.90}
\definecolor{Quartb}{rgb}{0.96,1,0.96}
\definecolor{Quinta}{rgb}{0.1,1,0.90}
\definecolor{Quintb}{rgb}{0.16,1,0.90}

\begin{document}


{\bf \Large 
{\begin{center}
Examples of  pairs of ordered congruent-like  $n$-gons  with different areas\\
\end{center}}
 }
 \bigskip
 {\begin{center}
 \large Michele Gaeta and Giovanni Vincenzi
 \end{center}
 \vskip0.5cm
 }
 

%


%
%
%
%
%
%
%
%


\address{Michele Gaeta \\
C.L. Matematica \\
Universit\`a di Salerno \\
via Giovanni Paolo II, 132 \\
I-84084 Fisciano (SA), Italy
  }
\email{m.gaeta23@studenti.unisa.it}

\address{Giovanni Vincenzi \\
Dipartimento di Matematica \\
Universit\`a di Salerno \\
via Giovanni Paolo II, 132 \\
I-84084 Fisciano (SA), Italy\\
MR Author ID: 329840 ,Orcid:  0000-0002-3869-885X
  }
\email{vincenzi@unisa.it}


{\color{brown}
}
\noindent Let's consider a cycle of $n$ rods of possible varying lengths ($n$ positive integer), attached with hinges. Clearly, if we have more than three rods   we can  move them, and  we can obtain different shapes ($SSS$ congruence doesn't hold past $n=3$). For example,  when $n=4$ we have a model, that an engineer would call \emph{articulated quadrilaterals}, that is  a `four-bar linkages'.  These have numerous practical applications. This is because if we fix the  side $AB$, then the paths which can be followed by $C$ and $D$ are clearly defined. 


\begin{figure}[htbp]
\centering
\includegraphics[width=0.6\textwidth]{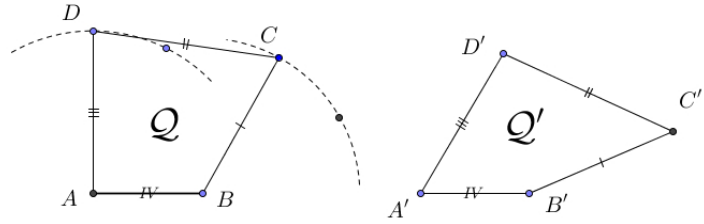} \quad
\caption{ Pair of articulated quadrilaterals. Their  sides are in the same order and have the same lengths; but  their inner angles are not pairwise congruent.}\label{fig:quadrilateri}
\end{figure}

\noindent What happens to a pair of articulated $n$-gons if we restrict the angles to be the same? 
Well, obviously if they are the same in the corresponding locations (so $SASASAS \ldots$ generalizes), then they will be congruent. 
What if we permute the angles? Here the answer is not so easy, because it depends on the number $n$ of the sides (angles). 

\begin{figure}[htbp]
\centering
\includegraphics[width=0.7\textwidth]{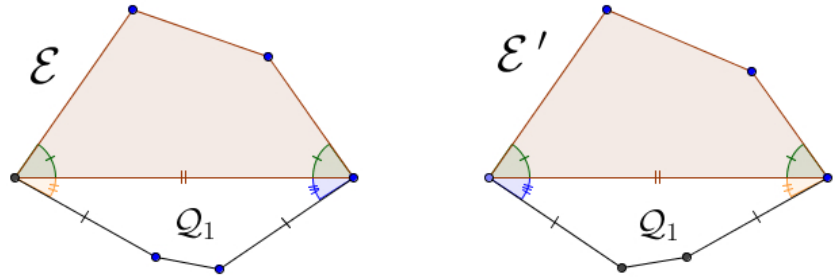}
\caption{
Pairs of ordered congruent-like hexagons $\mathcal E$ and $\mathcal E'$. They have been  constructed considering  two quadrilaterals, with a common side and suitable angles (green angles are congruent; orange angles are not congruent to 
blue angles). In the transition from left to  right, the quadrilateral $Q_1$ has been reflected, so $\mathcal E$ and $\mathcal E'$ are not congruent. 
}
\label{fig:esagoni}
\end{figure}

\begin{figure}[htbp]
\centering
\includegraphics[width=0.7\textwidth]{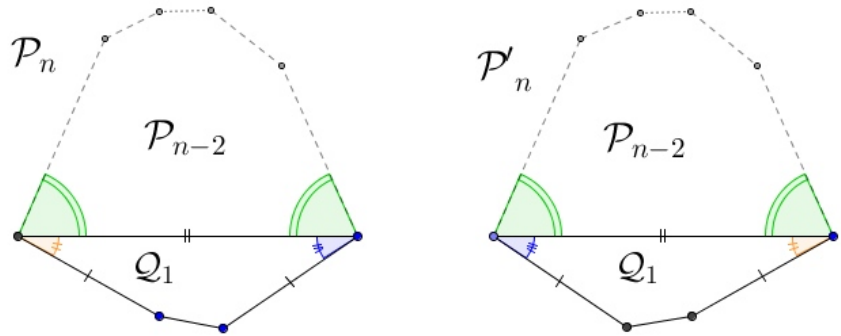}
\caption{
Pairs of ordered congruent-like $n$-gons $\mathcal P_n$ and $\mathcal P'_n$. 
}
\label{fig:n-gons}
\end{figure}

\smallskip
This kind of problems have been recently studied in  \cite{LaVi}, where the authors studied convex $n$-gons and introduced the formal definition of  
\emph{ordered congruent-like $n$-gons}  (see \cite[ p.47]{LaVi}). Here we highlight that the term ``ordered''  just refers to the sequence of the sides of the $n$-gons and not to their angles.  Intuitively, two ordered congruent-like $n$-gons can be considered as an  articulated  $n$-gon such that by moving its sides,   it is possible to find a new configuration  in which  every angle may change  position but  `preserves' its magnitude 
 (see  \cite[Figure 3]{LaVi}; 
  Figure  \ref{fig:esagoni} below or Figure  \ref{fig:n-gons}).



{\color{yellow}

}

 {\color{black}
 \noindent In \cite{LaVi} the authors   showed,  by an articulated proof which
splits into  many different sub cases corresponding to the possible permutations  of the angles,  
  that two  ordered congruent-like quadrilaterals are congruent. 
 The case $n=5$ is still open; while 
 even if it seems  strange at first glance, when $n>5$ it is not difficult to find 
 pairs of non-congruent $n$-gons  that are ordered congruent-like (see Figure \ref{fig:n-gons}). }
 
We note that the concept  of ``pairs of ordered congruent-like $n$-gons'' can be considered as a particular case of a more  general notion of
``pairs of  congruent-like $n$-gons'', in which both the order of the side  and the angles may be not preserved. 
Examples of ``pairs of  (non-ordered) congruent like $n$-gons'' that are not  congruent  can be  constructed  for every $n>3$ (see \cite[Figure 1]{AnLaVi} for an easy construction when  $n>4$). We   refer the reader  to  \cite {AnLaVi} and \cite{AnatrielloLaVi2019} for a deeper discussion of   (non-ordered) congruent-like polygons, and the case  $n=4$.

We highlight that at the moment all constructions of pairs of ordered  congruent-like $n$-gons' that  are known  produce pairs of polygons with equal areas. It is natural to ask: is this forced? 
We can summarize what we have  just resumed about ordered congruent-like $n$-gons  issues in the following scheme:


\begin{figure}[htbp]
\centering
\includegraphics[width=1.0\textwidth]{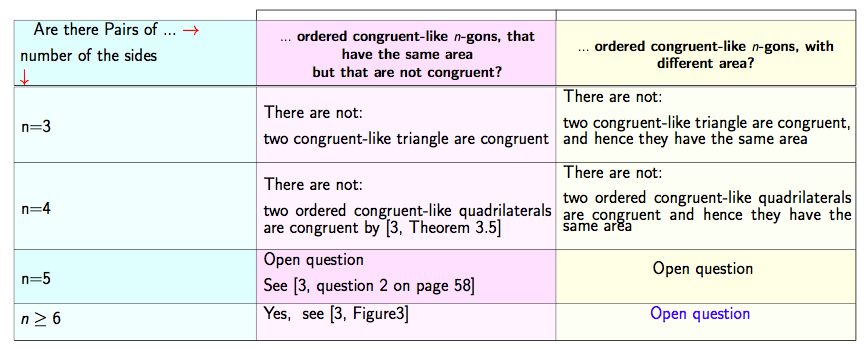} 
\caption{Summary related to pairs of ordered-congruent-like $ n $-gons.}\label{fig:step3}
\end{figure}

Now, we will see that  \emph{ for every $n>5$ there are pairs   of ordered congruent-like  $n$-gons  with different areas,} 
thus, the last cell of the above table may be fill with  ``Yes''.

\vskip1cm

\begin{figure}[htbp]
\centering
\includegraphics[width=0.6\textwidth]{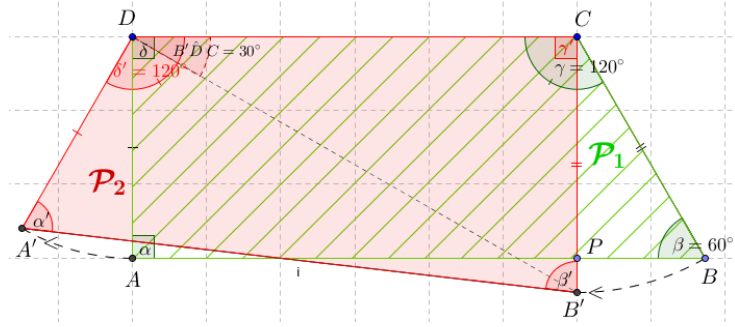}
\caption{ Quadrilateral $\mathcal P_2$ has been obtained from the trapezoid $\mathcal P_1$ by a rotation of $30^\circ$ of  the sides $AD$ and and $CB$.
}
\label{fig:step1}
\end{figure}

We begin the construction with the following trapezoid (see Figure \ref{fig:step1}).

  Starting from the  trapezoid $\mathcal P_1=(ABCD)$ that has been  obtained as a union of a rectangle
  with the sides of length 3 and 6
    and half equilateral triangle for which the height is 3,
     we can turn simultaneity  (clockwise) both the sides $DA$ and $CB$  around $D$ and $C$ (respectively) of $30^\circ$. Thus,  a new quadrilateral  $\mathcal P_2=(A'B'CD)$ is obtained.

\medskip
 First, we show that $Area (ABCD)< Area (A'B'CD)$.
 By   construction  $ CD=6$, $CB'=CB= 2 \sqrt 3$ and $D\hat C B'$ is right, so that 
 
$$CB'= CD \tan( B'\hat DC)\; {\mbox{and}}\;  DB'= \sqrt {(36 + 12)}=\sqrt {48}.$$

It follows that 
$ \tan( B'\hat DC) = \frac{2\sqrt 3}{6}=
\frac{1}{\sqrt 3}$, and hence $B'\hat DC= 30^\circ$. Therefore 
$ A'DB'$ is a right  triangle.
Now we can easily compute the areas of the two quadrilaterals:
\begin{itemize}
\item  $Area (ABCD)= APDC + \frac{1}{2}3 PB= 18+  \frac{1}{2}3 \sqrt 3$

\item $Area (A'B'CD)= CDH+ A'DH= \frac{1}{2}6\cdot  2\sqrt 3+\frac{1}{2}3 \sqrt {48}
= 6\sqrt 3+ 6\sqrt 3=12 \sqrt 3$.
\end{itemize}
It follows that $Area (ABCD)<Area (A'B'CD)$.

We also note that $\alpha'+\beta'=150^\circ$.
 
\begin{figure}[htbp]
\centering
\includegraphics[scale=0.30]{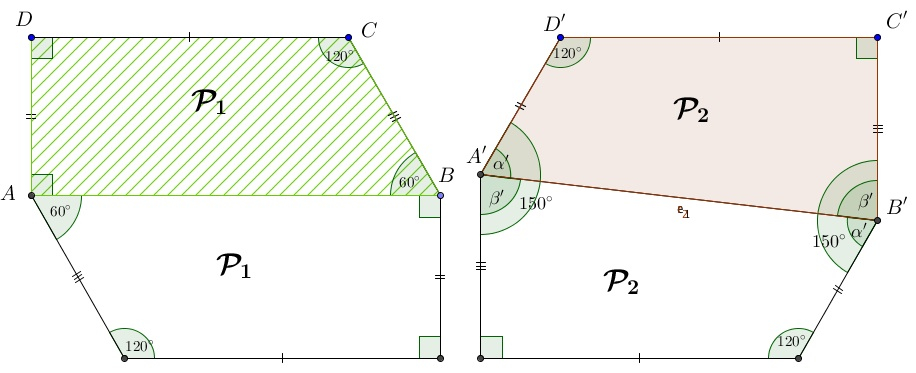}
\caption{Pairs of ordered congruent-like hexagons with different area. }\label{fig:step2}
\end{figure}

\medskip
 Now, to determine a pair of  hexagons as required, it is enough to consider the union of two copies of $\mathcal P_1$ and two copies of $\mathcal P_2$ as shown in Figure \ref{fig:step2}.
%
 

\medskip
\medskip

 The general case ($n>6$) can be easily derived from the hexagons adjoining  a suitable ($n-4$)-gon (see $\mathcal P_3$ in  Figure \ref{fig:step3}). 

 
 

\begin{figure}[htbp]
\centering
\includegraphics[scale=0.30]{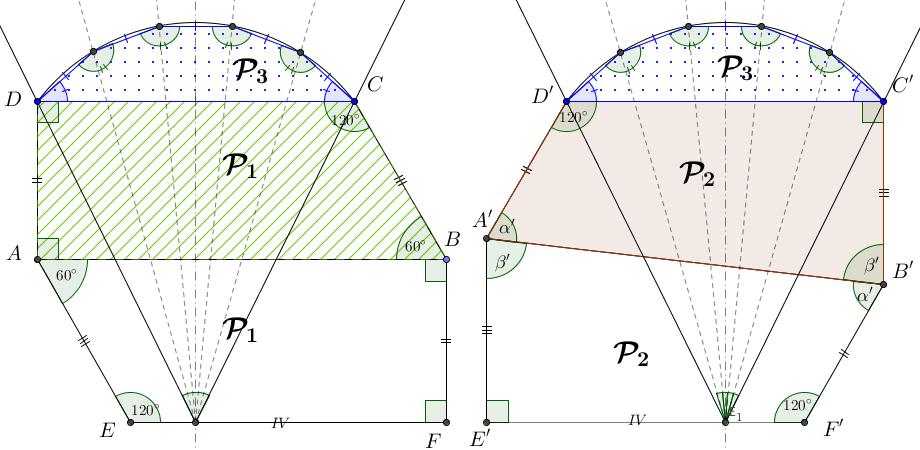}
\caption{Pairs of ordered congruent-like $n$-gons  with different areas.}\label{fig:step3}
\end{figure} 
\medskip
\medskip

As an  exercise,  the reader  can try to find two congruent-like heptagons with different areas.  


\vskip0.5cm

\begin{remark} It should be noted that the  polygons  constructed as above, say $\mathcal V_n$ e $\mathcal V'_n$,  are very special pairs of congruent-like polygons: they are  articulated $n$-gons  with different areas  with the property  that they can be obtained  from each  other by a deformation that preserves the order and the length of their  sides, and such that {each ordered (clockwise) sequence} of the consecutive  inner angles of $\mathcal V_n$ is equal to  some ordered (clockwise or counterclockwise) sequence of the consecutive  inner angles of $\mathcal V'_n$. 
{\color{black} Pairs of polygons of this kind can provide some useful hints  for a  lesson on  polygons. For example, it could be observed that  to determine the area  of a $n$-gon, it is not enough to know all its  sides and  angles; but we also need to know which are their  mutual positions. The use of dynamic geometry software could be useful in this study.}
\end{remark}

We conclude by observing that it does not seem easy to find similar constructions for pentagons. 


\vskip1cm
{\color{black}
We wish to thank the anonymous referees for their suggestions and comments,
which have improved the paper with respect  to the first submission.
}

\bibliographystyle{plain}

\end{document}